\documentclass[a4paper]{IEEEtran}
\usepackage{amsfonts,amssymb}
\usepackage{graphicx,epstopdf}
\usepackage{tikz}
\usetikzlibrary{automata}
\usepackage{color}
\usepackage{threeparttable}
\usepackage{amsmath}
\usepackage{amsfonts,amssymb}
\usepackage{slashbox} 
\usepackage{epstopdf}

\newtheorem{proposition}{Proposition}

\usepackage{algorithm}
\allowdisplaybreaks[4]

\IEEEoverridecommandlockouts
\begin{document}
\title{Modelling and Traffic Signal Control of Heterogeneous Traffic Systems}

\author{
%	Yicheng Zhang, Rong Su, Kaizhou Gao, Chunyang Sun and Christos Cassandras% <-this % stops a space
	Yicheng~Zhang,~\IEEEmembership{Student~Member,~IEEE,} Rong~Su,~\IEEEmembership{Senior~Member,~IEEE,} Yi~Zhang,~\IEEEmembership{Student~Member,~IEEE,} and Chunyang~Sun
	%%%%%%%%%%%%%%%%%%%%%%%%%%%%%%%%%%%%%%%%%%%%%
	\thanks{This research is supported by the research grand S15-1105-RF-LLF URBAN
		from the Economic Development Board, Singapore, for the project of
		``Development of NTU/NXP Smart Mobility Test-bed".}% <-this % stops a space
	\thanks{$^{1}$Yicheng Zhang, Rong Su Yi Zhang and Chunyang Sun are affiliated with School of Electrical and Electronic Engineering at Nanyang Technological University, Singapore 639798. E-mails
		{\tt\small yzhang088@e.ntu.edu.sg, rsu@ntu.edu.sg, yzhang120@e.ntu.edu.sg, suncy@ntu.edu.sg}}%
%	\thanks{$^{2}$Christos Cassandras is affiliated with the Division of Systems Engineering at Boston University, Brookline, MA 02446. E-mail:
%		{\tt\small cgc@bu.edu}}    %
}

%\author{Yicheng Zhang}
%\author
%{
%	\IEEEauthorblockN{Yicheng Zhang}
%	\IEEEauthorblockA
%	{
%		School of Electrical and Electronic Engineering\\
%		Nanayang Technological University, Singapore\\
%		yzhang088@e.ntu.edu.sg
%	}
%}
%% make the title area
\maketitle

\begin{abstract}
	An urban traffic system is a heterogeneous system, which consists of different types of intersections and dynamics. In this paper, we focus on one type of heterogeneous traffic network, which consists of signalized junctions and non-signalized ones, where in the latter case vehicles usually follow the first-in-first-out principle. We propose a novel model describing the dynamic behaviors of such a system and validate it via simulations in VISSIM. Upon such a new model, a signal control problem for a heterogeneous traffic network is formulated as a mixed integer programming problem, which is solved by a Lagrangian multiplier based hierarchical distributed approach. Comparisons between a homogeneous traffic system and a heterogeneous one are provided, which leaves the door open for developing a systematic planning approach on deciding what traffic junctions require signal control to ensure a good traffic control performance, thus, have a great social and economic potentials, considering that it is rather expensive to have signal control in an urban area.
%The urban traffic system is a heterogeneous system consisting of multiple kinds of intersections and dynamics. In this paper, a novel model for the urban traffic system is proposed to describe the system with both signalized intersections and non-signalized intersections. The proposed model for the non-signalized intersections with uncontrollable flows is validated based on the simulations in VISSIM. Moreover, a signal control problem formulation for the traffic network with signalized and non-signalized intersections is proposed as a mixed integer programming problem. Some comparisons between the fully controllable traffic system and partially controllable traffic system are provided and we show some potential applications for the traffic system design of this methodology.
\end{abstract}

%\IEEEpeerreviewmaketitle
\section{Introduction}

Traffic congestions in urban area occur frequently, which affect daily life and pose all kinds of problems and challenges. Alleviation of traffic congestions not only improves travel safety and efficiencies but also reduces environmental pollution.
The urban traffic system consists of intersections and links, which include controllers, flow dynamics and volume dynamics. 
The traffic system modelling and the traffic signal design are the two important phases in solving the urban traffic congestion problem.

From the system modelling side, the urban traffic system is a heterogeneous system consisting of multiple kinds of intersections, dynamics and participants. 
%In this paper, different kinds of intersections and the traffic flow dynamics at these intersections are taken into consideration. 
In general, the intersections in urban traffic system could be classified into two main groups, i.e., the signalized intersections and the non-signalized intersections. A signalized intersection equips a set of traffic signals while 
a non-signalized or uncontrolled intersection is one in which the entrance into the intersection from any of the approaches is not controlled by a regulatory sign or a traffic signal.
Although it seems that putting traffic signals at all intersections in the urban area would lead to a very orderly fashion, the costs of building and maintaining signalized intersections are much higher than the costs of non-signalized intersections. 
Moreover, it may be unnecessary to signalize all the intersections because of the low traffic flow for some minor intersections. 
Thus in real-world practice, both signalized and non-signalized intersections are adopted in the urban traffic system, which lead to the heterogeneous features of the link volume dynamics and the traffic flow dynamics.
A detailed survey on monitoring an intersection based on the bahaviors of the participants is shown in \cite{Shirazi2017}. It provides a microscopic view or Lagrangian Model of the traffic dynamics at the intersection.
For a macroscopic view or Eulerian Model of the traffic system, some descriptions about the non-signalized intersections or uncontrolled intersections can be found in the literature \cite{Rob69}-\cite{WuQS2005}.

From the traffic signal scheduling strategy design view, there are four different traffic signal control strategies,  namely, fixed  time  strategies  versus  traffic  responsive  strategies,  and  isolated  strategies versus coordinated strategies. Notable strategies proposed in the last few decades include, SIGSET \cite{All71}, SIGCAP\cite{All76}, TRANSYT \cite{Rob69} \cite{LG99}, SCOOT \cite{HRB82} \cite{Rob91}, OPAC \cite{Gar83} and PRODYN \cite{FHT83}. From the system modelling view, Lighthill and Whitham in 1955 and Richards in 1956 provided a good foundation in macroscopic traffic system modelling, the LWR model \cite{LWR55} \cite{LWR56}.
A breakthrough came when Daganzo developed a finite difference solution scheme for the LWR model by adopting a simplified Fundamental Diagram \cite{Dag94} \cite{Dag95}, which he called the cell-transmission model (CTM). 
In \cite{lo1999novel} and \cite{lo2001cell}, the author transforms CTM into a mixed integer linear programming (MILP) problem. However, the transformation is not equivalent rigorously. Moreover, the formulation does not include any analysis about the relationship among the speed, previous signal status and velocity. Another weak point is that in this model the cycle time is fixed and the signal scheduling strategy is based on the fixed-time cycle. 
Lin and Wang \cite{lin2004enhanced} define difference cell models in the traffic network and different constraint sets are proposed for each type of cells.
In literature \cite{Zhang2015}, an urban road traffic light control problem is formulated as a scheduling problem, aiming to reduce the total waiting time over a given finite horizon. One of the key contributions in the model is to describe each outgoing flow rate as a nonlinear mixed logical switching function over the source link volume, the destination link volume and capacity, and the driver’s potential psychological response to the past traffic light signals. The outgoing flow rate model makes the proposed approach applicable to both under-saturated and over-saturated situations. The traditional concepts of cycles, splits, and offsets are not adopted in this framework, making the proposed approach fall in the class of model-based optimization methods, where each traffic light is assigned with a green light period in a real-time manner by the network controller.  
Previously the traffic signal control (TSC) problem is solved by using fixed time strategy with individual or coordinated intersections to adjust the green signal ratio in the traffic system.
It can be formulated as an optimization problem and solved by standard optimization methodologies or computational intelligence techniques \cite{Zhao2012}.

In this paper, a novel model for the urban traffic system is proposed to describe the heterogeneous traffic system with both signalized intersections and non-signalized intersections. For the signalized intersections, the proposed model describes the flow dynamics for controllable traffic flows and uncontrollable traffic flows, e.g., some merging flows from the free-turning directions. For the non-signalized intersections, the model provides a method to estimate the flow dynamics based on the gap-acceptance model and first-come-first-serve principle. The proposed non-signalized intersection model is validated based on the simulations in VISSIM simulator, a microscopic traffic simulation platform \cite{VISSIM2010}. 
Moreover, a signal control problem formulation for the heterogeneous traffic network is proposed as a mixed integer programming problem. 
In the simulation section, some comparisons among the fully controlled traffic system (all intersections are signalized), partially controlled traffic system and fully uncontrolled traffic system (all intersections are non-signalized) are provided, which leaves the door open for developing a systematic planning approach on deciding what traffic junctions require signal control to ensure a good traffic control performance.

The contributions of this work are listed below. Firstly, we develop and validate a heterogeneous traffic network model with both signalized and non-signalized intersections; secondly, a traffic signal control problem is formulated based on this model to minimize the network-wise traffic delay; thirdly, we do some simulations to to provide some guidance for the design of urban traffic system.

The paper is organized as follows. In Section II, we introduce the urban traffic system model for the problem development. Mathematical models for both signalized intersections and non-signalized intersections are proposed in this section and all the constraints for the flow dynamics and volume dynamics are formulated. An urban traffic signal control (UTSC) problem formulation for this heterogeneous system model is provided in Section III and a distributed way to solve this UTSC problem based on Lagrangian multiplier method is shown in Section IV. The simulation-based model validation as long as the simulation results from both numerical analyses in MATLAB \cite{Matlab2015} with Gurobi \cite{Gurobi2016} and simulations in VISSIM are provided in Section V. Some conclusions are drawn in Section VI.
 
\section{Modelling of Traffic Network with Signalized and Non-Signalized Intersections}
An urban traffic network consists of a set of road \textit{links} connecting with each other via \textit{intersections}. Each intersection consists of a number of approaches and the crossing area. An approach may have one or more lanes but has a unique, independent queue. Approaches are used by corresponding traffic streams. Two \textit{compatible} streams can safely cross the intersection simultaneously, while \textit{antagonistic} streams cannot. In traditional traffic signal control, a signal cycle is one repetition of the basic series of \textit{stages} at an intersection, where each stage consists of simultaneous traffic light signals allowing a predefined compatible traffic streams to cross the intersection simultaneously. The duration of a cycle is called \textit{cycle time}. For safety reasons, constant lost (or intergreen) times of a few seconds are necessarily inserted between consecutive stages to avoid interference between antagonistic streams. For each traffic light, the ratio of the green and red times within one cycle is called \textit{split}, and a delay between the starting times of green periods of two neighboring traffic lights along the same traffic route is called \textit{offset}.

\subsection{A traffic network model}\label{subsectiontraffic}
\begin{figure}[!htb]
	\centering
	\includegraphics[width=\linewidth]{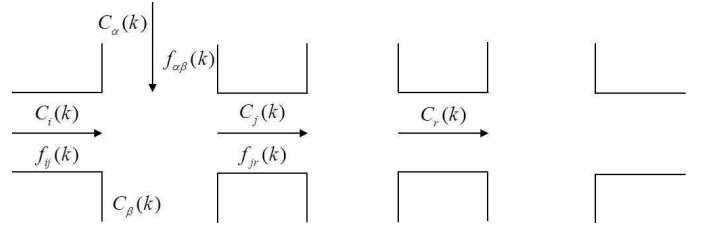}
	\caption{A traffic network}
	\label{fig01}
\end{figure}

A traffic network consists of a set of road links and intersections, i.e., both signalized intersections and non-signalized intersections. Fig. \ref{fig01} depicts a simple unidirectional traffic network, where each intersection has only two antagonistic traffic streams (or flows).
We adopt a discrete-time model similar to a cell transmission model \cite{Dag94}. Let $\Delta$ be a given sampling interval, e.g., $\Delta = 5s$. We use
$C_i(k)$ to denote the {\color{black} integer} number of vehicles (or the \emph{volume}) in the link $i$ during the time interval $k$, and $f_{ij}(k)\geq 0$ for the {\color{black} integer} number of vehicles (or the \emph{shift}) from the link $i$ to the link $j$ in the interval $k$. Each link is partitioned into \emph{segments}, where each vehicle can (approximately) traverse a segment in $\Delta$ with a free flow speed. 
%For example, in Singapore the speed limit of all urban roads is 50km/h (unless otherwise stated), thus, the length of a segment within $5s$ is 69.4m. Considering that a driver may drive at about  35km/h when starting to move or making a turn in a free flow, the length of a segment can be between about 50m and 70m. Let $\mathbb{N}$ denote the set of natural numbers. We have the following result.
%\begin{proposition}\label{prop1}
%	\textnormal{For any natural number $n\geq \mathbb{N}$ no smaller than 10, there exists natural numbers $k_1,k_2,k_3\in\mathbb{N}$ such that $n=5k_1+6k_2+7k_3$.\hfill $\Box$}
%\end{proposition} 
%Proof: For any natural number between 10 and $19$ we have:
%\begin{eqnarray*}
%	\begin{aligned}
%		& 10 = 5+5,\,\, 11=5+6,\,\, 12=5+7,\,\,13=6+7,\\
%		& 14=7+7,\,\, 15 = 5+5+5,\,\, 16=5+5+6, \\
%		& 17=5+5+7,\,\, 18=5+6+7,\,\, 19=5+7+7.
%	\end{aligned}
%\end{eqnarray*}   
%For any $n'\in\mathbb{N}$ no smaller than 20, we can represent it as $n'=n+10k$, where $10\leq n\leq 19$ and $k\in\mathbb{N}$. Thus, the propagation follows.\hfill $\blacksquare$

%If we round up the length of each link to its second digit, e.g., to represent 101m as 100m, and 105m as 110m, etc., then by Prop. \ref{prop1} we know that each link longer than 100m can be represented by a few segments, whose lengths are either 50m, 60m  or 70m. We do not consider a link, whose length is below 100m, as this is usually quite rare in practice.   

Let $\mathcal{L}$ be the set of all such one-way segments obtained above,   and $\mathcal{J}_r$ the set of all (real) intersections. For each pair of consecutive segments belong to the same link, we consider them to be connected by a virtual intersection, where the corresponding traffic light is alway GREEN. We denote the set of all those virtual intersections as $\mathcal{J}_v$ and let $\mathcal{J}:=\mathcal{J}_r\cup\mathcal{J}_v$. With a slight abuse of terminology, from now on we still call each segment a link. Two groups of intersections are involved in the system, i.e., the signalized intersections and non-signalized intersections. We use $\mathcal{J}_s$ and $\mathcal{J}_n$ to denote the signalized and non-signalized ones, respectively. For each signalized intersection $J\in\mathcal{J}_s$, let $\Omega_J$ be the set of stages in the intersection $J$, and $\mathcal{F}_J\subseteq \mathcal{L}\times \mathcal{L}$ be the set of all streams in the intersection $J$, i.e., $(i,j)\in \mathcal{F}_J$ means that there exists a traffic stream from the link $i$ to the link $j$ via the intersection $J$. Let $h_J:\Omega_J \rightarrow 2^{\mathcal{F}_J}$ be the association of each stage to relevant compatible streams. Here, for simplicity we assume that for any two different stages $w_i$ and $w_j$,  $h_J(w_i)\cap h_J(w_j)=\varnothing$. For the non-signalized intersection $J\in\mathcal{J}_n$, $\mathcal{F}_J\subseteq \mathcal{L}\times \mathcal{L}$ is the set of all streams in the intersection $J$ and there is no stage for such kind of intersections.
For the traffic dynamics in the urban traffic system, we make the following assumption, which is suitable for a deterministic analysis.
\begin{itemize}
	\item \textbf{A1:} The traffic demand (i.e., vehicle entrance and exit) models are known, and there is an infinite queue for each entrance demand.
	\item \textbf{A2:} The link turning ratios in the network are known.
	\item \textbf{A3:} Each vehicle inside the network will leave the network, delayed only by traffic signals.
\end{itemize}
Assumption \textbf{A1} specifies that the number of vehicles requiring to enter or leave the network at relevant locations in the network is known, mainly due to historical data adjusted by real-time model identification, which is out of the scope of this paper. Once an entrance demand appears, it will not disappear until it is served, i.e., there is an infinite queue at each demand location to hold all vehicles, which require to enter the network. Assumption \textbf{A2} can be ensured by using historical data adjusted by real-time data. Assumption \textbf{A3} simply says that no vehicle will stop in the network unnecessarily, i.e., whenever allowed by the traffic signals, it will move forward. This assumption is reasonable for all normal traffic situations.  
In the following sections, we provide the models for the links, signalized intersections and non-signalized intersections, respectively.   

\subsection{Link dynamics constraints}
Due to the conservation of vehicles, each link $j\in\mathcal{L}$ has the following volume dynamics:
\begin{equation}
(\forall k\in\mathbb{N})\, C_j(k+1) = C_j(k) + d_j(k)-s_j(k),
\label{eqn_vehicle_1}
\end{equation}
{\color{black}where $d_j(k)$ and $s_j(k)$ are the number of entrance vehicles and exit vehicles of the link $j$ at $k$ respectively, i.e.,}
\begin{subequations}
	\begin{align}
	d_j(k)&=\sum_{i\in\mathcal{L}:(i,j)\in\cup_{J\in\mathcal{J}}\mathcal{F}_J}f_{ij}(k)\\
	s_j(k)&=\sum_{i\in\mathcal{L}:(j,i)\in\cup_{J\in\mathcal{J}}\mathcal{F}_J}f_{ji}(k)
	\end{align}
\end{subequations}
For the example shown in Fig. \ref{fig01}, $d_j(k)=f_{ij}(k)$ and $s_j(k)= f_{jr}(k)$. 
We define the capacity of each link based on the length of link and the minimal separations between vehicles. For example, the length of the link is 200m, the average vehicle length is 5m and the minimal separation between each two vehicles is 1m, then the capacity is about $\lfloor200/(5+1)\rfloor = 33$ (veh). Denote the capacity for link $j$ as $\hat{C}_j$,
\begin{equation}
(\forall k\in\mathbb{N})\, C_j(k) \leq \hat{C}_j,
\end{equation}
Note that no matter the link is linked to the signalized or non-signalized intersection, the equation for the link volume dynamics will always hold. Thus we do not require to distinguish the characteristics for the intersections connecting to the link. However, the traffic flow dynamics to describe the vehicle's approaching or departing from the signalized or non-signalized intersections are different and we propose different traffic flow dynamics for them in our model.

\subsection{A Model for Signalized Intersections}
\subsubsection{Stage constraints}
In each given interval $k$, there exists only one active stage for an intersection $J\in\mathcal{J}_r$, which is captured by the following constraints:
\begin{subequations}
	\begin{align}
	& (\forall w\in\Omega_J)\, \theta_w(k)=0\Rightarrow (\forall (i,j)\in h_J(w))\, f_{ij}(k)=0
	\label{eqn_RG_1}\\
	&  \sum_{w\in\Omega_J}\,\theta_w(k)=1
	\label{eqn_RG_2}\\
	& (\forall w\in\Omega_J)(\forall k\in\mathbb{N})\, \theta_w(k)\in\{0,1\}
	\label{eqn_RG_3}
	\end{align}
\end{subequations}
where $\theta_w(k)=0$ and $\theta_w(k)=1$ denote the RED and GREEN traffic lights associated with the stage $w$  respectively. Condition (\ref{eqn_RG_1}) indicates that if the stage traffic light is RED, then all relevant shifts are zero. Condition (\ref{eqn_RG_2}) indicates that there can be only one GREEN traffic stage at any time. For all $J\in\mathcal{J}_v$ and all $w\in\Omega_J$ we have $\theta_w(k)=1$.

%If $l_{ij}(k)$ is determined, the link outgoing shift $f_{ij}(k)$ is given as follows:
%{\color{black}\begin{equation}
%f_{ij}(k) = \lfloor \min\{\lambda_{ij}(k)C_i(k),l_{ij}(k)v_i^*d^*\Delta,\hat{C}_j-C_j(k) \}\rfloor
%		\label{eqn_fr_1}
%\end{equation}}
%where $\lfloor\cdot\rfloor$ is the largest integer not greater than input argument, $\lambda_{ij}(k)$ is the turning ratio of vehicles in the link $i$ towards the link $j$ at $k$, which is assumed to be known in advance, 
%{\color{black}$v_i^*$ is the free speed in the link $i$, and $d^*$ is the maximum link density, i.e., the density of the situation where all vehicles are considered having the same standard length with the minimum separation distance.}   
%Clearly, $\sum_{j\in\mathcal{L}:(i,j)\in\cup_{J\in\mathcal{J}}\mathcal{F}_J}\lambda_{ij}(k)=1,\label{eqn_turing_ratio}$ meaning that each vehicle in the link $i$ will move into some downstream link $j$. {\color{black}Conditions (\ref{eqn_fr_1}) indicate that the number of vehicles in one time interval, $f_{ij}(k)$, is the largest integer that is upper bounded by the upstream volume $\lambda_{ij}(k)C_i(k)$ of the link $i$, the downstream remaining capacity $\hat{C}_j-C_j(k)$, and the expected shift attainable by discounting the maximum shift $f_{i,max}:=v_i^*d^*\Delta$ by the speed level $l_{ij}(k)$. }
\subsubsection{Flow dynamics constraints}
For each stage $w\in \Omega_J$ and each stream $(i,j)\in h_J(w)$, the exit shift $f_{ij}(k)$ is determined by the current upstream link volume $C_i(k)$, the current remaining downstream link capacity $\hat{C}_j-C_j(k)$, where $\hat{C}_j$ is the capacity of link $j$, and the traffic light signals in the past $r+1$ time intervals $\theta_w(k-r),\cdots,\theta_w(k)$, i.e.,
\begin{equation}\label{fijlaw01}
f_{ij}(k) = g_{ij}(C_i(k),\hat C_j-C_j(k), \theta_w(k-r), \dots, \theta_w(k)),
\end{equation}
where $g_{ij}(\cdot)$ is a nonlinear function. The motivation behind this model is that if the stage $w$ has been active for the past $r$ intervals, then the drivers intend to keep a high speed $v_{ij}(k)$ as long as the downstream link has sufficient capacity to receive the flow, i.e., $v_{ij}(k)=l_{ij}(k)v_i^*$, where the coefficient $l_{ij}(k)\in [0,1]$ is monotonically nondecreasing with respect to the length of a GREEN period, except for the last few moments (or an amber or yellow light period in the current practice) when the driver needs to slow down to a complete stop to anticipate the beginning of a RED period, and $v_i^*$ is the maximum free flow speed. 

%As an illustration, Figure \ref{fig0} depicts the driving speed observed at the intersection of General Holms Drive and Bestic Street in Sydney, Australia \cite{Akcelik2008}.
%
%\begin{figure}[!tbh]
%	\centering
%	\includegraphics[width = \linewidth]{Speed-profile.png}
%	\caption{Observed driving speed profile at one intersection when a green light is on}
%	\label{fig0}
%\end{figure}

\begin{figure}[tbh]
	\centering
	\includegraphics[width=0.8\linewidth]{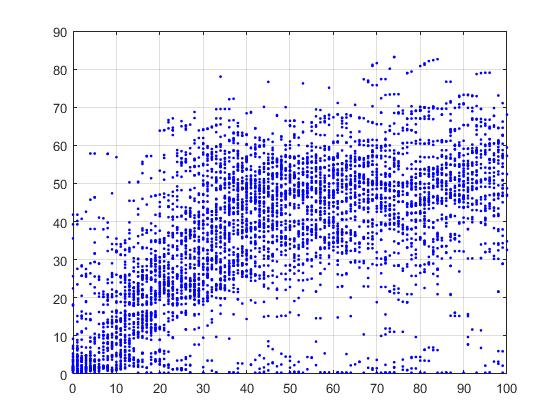}
	\includegraphics[width=0.8\linewidth]{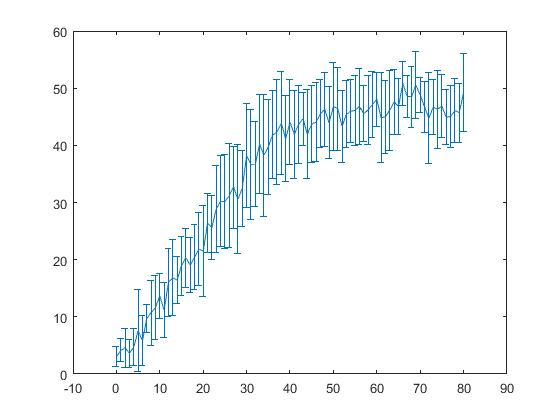}
	\includegraphics[width=0.8\linewidth]{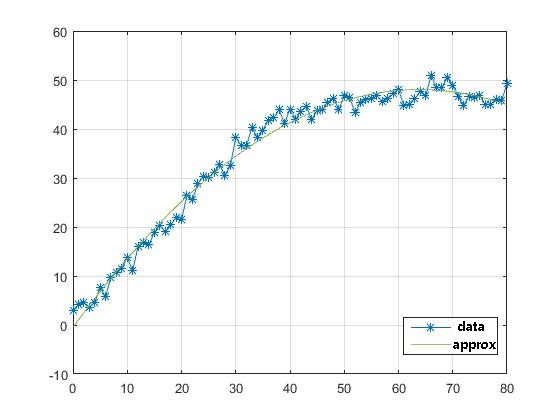}
	\caption{Relationship between the speed and the green time length shown in road test}
	\label{fig:1011amt12}
\end{figure}

We have also done some road tests in Singapore with speed guns and the results are shown in Fig. \ref{fig:1011amt12}.
The data is collected at an intersection in Singapore by two speed guns. We select four time intervals during two days which include the peak hours and off-peak hours and eight groups of data are recorded.
The leftmost picture shows the raw data, which includes about 30000 points. The middle picture shows the processed data with error bars and the  rightmost picture shows the processed data with curve fitting.
In the middle and right pictures, the x-axis shows the green time length at the intersection and the y-axis shows the average speed for vehicles passing this intersection. We can conclude that with the increasing of the length of green time at the intersection, the average exit speed is increasing and finally saturated with a certain value which is always determined by the upper speed limit.

It is clear that the average exit speed in consecutive time intervals can be approximated by a set of discrete speed levels, say $l_{ij}^0\geq \cdots\geq l_{ij}^r>0$. Due to the monotonically nondecreasing nature of $l_{ij}(k)$ with respect to the length of a GREEN period mentioned before (except for the last GREEN interval before a RED period starts), we introduce the following generic way of defining the speed level $l_{ij}(k)$ based on the expected traffic light assignment:

\begin{subequations}
	\begin{align}
	& \theta_w(k+1)=1\Rightarrow  l_{ij}(k)=\sum_{p=0}^r\delta_{ij}^p(k)l_{ij}^p
	\label{eqn_speed_1}\\
	&  \theta_w(k+1)=0
	\Rightarrow    l_{ij}(k)=\frac{1}{2}l_{ij}(k-1)
	\label{eqn_speed_1.5}\\
	& \sum_{p=0}^r\delta_{ij}^p(k)-\theta_w(k) =0
	\label{eqn_speed_2}\\
	\nonumber&  (\forall q:0\leq q\leq r)\,(1-\theta_w(k-q-1))\\
	& \prod_{p=0}^q\theta_w(k-p)=1\wedge \theta_w(k+1)=1 
	\Leftrightarrow \ \   \delta_{ij}^{r-q}(k)=1
	\label{eqn_speed_3}\\
	& (\forall p:0\leq p\leq r)\, \delta_{ij}^p\in\{0,1\}
	\label{eqn_speed_5}
	\end{align}
\end{subequations}

Condition (\ref{eqn_speed_1}) says that, as long as $k$ is not the last GREEN interval, i.e., $\theta_w(k+1)=1$, the corresponding speed level $l_{ij}(k)$ will be determined by the subsequent conditions. If the next time interval is RED, i.e., $\theta_w(k+1)=0$, as stated by Condition (\ref{eqn_speed_1.5}), then the speed level is one half of the speed level in $k-1$, assuming that the traffic speed needs to be reduced linearly from $l_{ij}(k-1)$ to 0 in $\Delta$.  Condition (\ref{eqn_speed_2}) indicates that if the traffic light of the stage $w$ is RED, i.e., $\theta_{w}(k)=0$, then $\sum_{p=0}^r\delta_{ij}^p(k)=0$, which by Condition (\ref{eqn_speed_1}) means $l_{ij}(k)=0$, i.e., the speed in the stage $w$ must be zero; if the traffic light of the stage $w$ is GREEN, i.e., $\theta_{w}(k)=1$, then by Condition (\ref{eqn_speed_1}), $l_{ij}(k)$ can only choose one speed level because of $\sum_{p=0}^r\delta_{ij}^p(k)=1$. Conditions (\ref{eqn_speed_3}) indicates that the actual speed level depends on the number of consecutive green light intervals from $k$ backward in time - the larger the number of consecutive green intervals, the higher the speed level. 

If $l_{ij}(k)$ is determined, the link outgoing shift $f_{ij}(k)$ is given as follows:
\begin{equation}
\begin{aligned}
& f_{ij}(k) = \lfloor \min\{\gamma_{ij}(k)C_i(k),\hat{C}_j-C_j(k)+s_j(k),\\
& \qquad \qquad l_{ij}(k)v_i^*d^*\Delta \}\rfloor
\end{aligned}
\label{eqn_fr_1}
\end{equation}
where $\lfloor\cdot\rfloor$ is the largest integer not greater than input argument, $\gamma_{ij}(k)$ is the turning ratio of vehicles in the link $i$ towards the link $j$ at $k$, which is assumed to be known in advance, 
$v_i^*$ is the free speed in the link $i$, and $d^*$ is the maximum link density, i.e., the density of the situation where all vehicles are considered having the same standard length with the minimum separation distance.
Clearly, $\sum_{j\in\mathcal{L}:(i,j)\in\cup_{J\in\mathcal{J}}\mathcal{F}_J}\gamma_{ij}(k)=1,\label{eqn_turing_ratio}$ meaning that each vehicle in the link $i$ will move into some downstream link $j$. Conditions (\ref{eqn_fr_1}) indicate that the number of vehicles in one time interval, $f_{ij}(k)$, is the largest integer that is upper bounded by the upstream volume $\gamma_{ij}(k)C_i(k)$ of the link $i$, the downstream remaining capacity plus the exiting vehicles which are current in the downstream link, i.e.,  $\hat{C}_j-C_j(k)+s_j(k)$, and the expected shift attainable by discounting the maximum shift $f_{i,max}:=v_i^*d^*\Delta$ by the speed level $l_{ij}(k)$.

\subsubsection{Merging constraints for signalized intersections}
\begin{figure}[!ht]
	\centering
	\includegraphics[width = 0.8\linewidth]{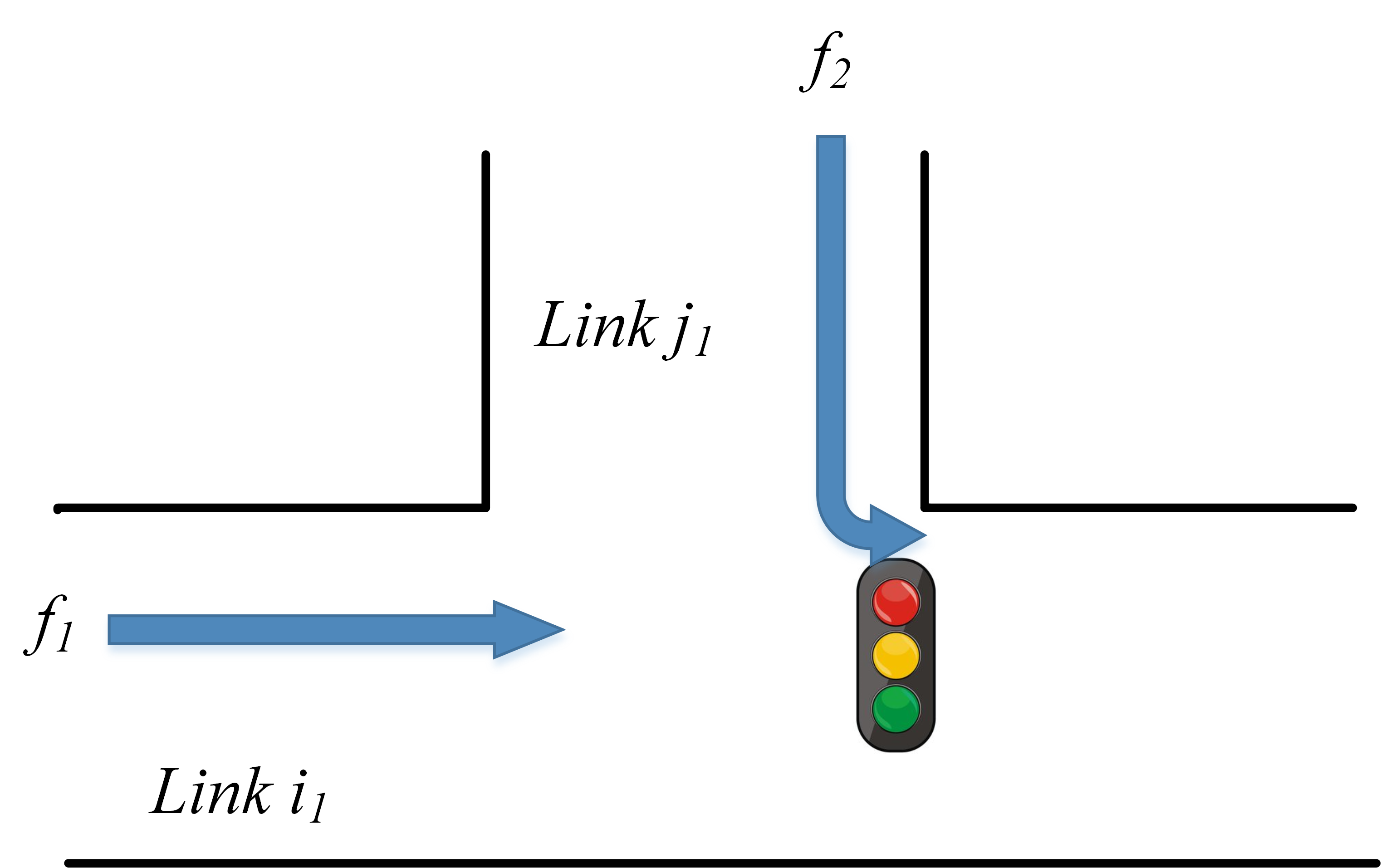}
	\caption{An example for traffic flow merging}
	\label{fig_merging}
\end{figure}

Mergence of multiple flows in one link is commonly seen in a traffic network. For example, Fig \ref{fig_merging} shows the mergence of traffic flows from left to right horizontally and from top to right (with a left turn).

In this case, based on the traffic rule \cite{TrafficRule2016}, a vehicle at an uncontrolled road intersection must \textit{give way} to any vehicle or stream of vehicles immediately approaching him from its right or offside. For each intersection $J\in\mathcal{J}_s\cup \mathcal{J}_n$ let $P_J:\mathcal{F}_J\rightarrow\mathbb{N}$ be a \emph{priority} map such that for any two streams $\omega_1,\omega_2\in\mathcal{F}_J$, $P(\omega_1)<P(\omega_2)$ implies that stream $\omega_1$ has a higher priority than that of $\omega_2$. With this interpretation, we have the following rule: 
\begin{equation}
\begin{aligned}
& (\forall (i,j),(i',j)\in \mathcal{F}_J)(\forall k\in\mathbb{N}) \\
&\qquad P(i,j)<P(i',j)\wedge f_{ij}(k)\neq 0 \Rightarrow f_{i'j}(k) = 0,
\end{aligned}
\label{eqn_merging}
\end{equation}
namely, when two streams $(i,j)$ and $(i',j)$ merge into the link $j$, if $(i,j)$ has a higher priority than $(i',j)$, then there exists a non-empty flow in $(i',j)$ only if there exists no flow in $(i,j)$, i.e., vehicles in $(i',j)$ needs to give way to vehicles in $(i,j)$.
The actual flow rate in $(i',j)$ during the merging process is described as follows:
\begin{equation}
f_{i'j}(k) \leq \hat{C}_j - \big[C_j(k) - s_j(k) + f_{ij}(k)\big]
\label{eqn_merging3}
\end{equation}
The inequality (\ref{eqn_merging3}) means that for ant time interval $k$, the give-way flow should be no more than the difference between the capacity of the downstream link $j$ and its current link volume, which takes the current incoming flow with right-of-way $f_{ij}(k)$ into account.

\subsection{A Model for Non-Signalized Intersections}
For each non-signalized intersection, we consider only the type of all-way stop-controlled intersection (AWSC), where each approaching traffic flow  is required to have a full stop before going through the intersection. 
\begin{figure}[!ht]
	\centering
	\includegraphics[width=\linewidth]{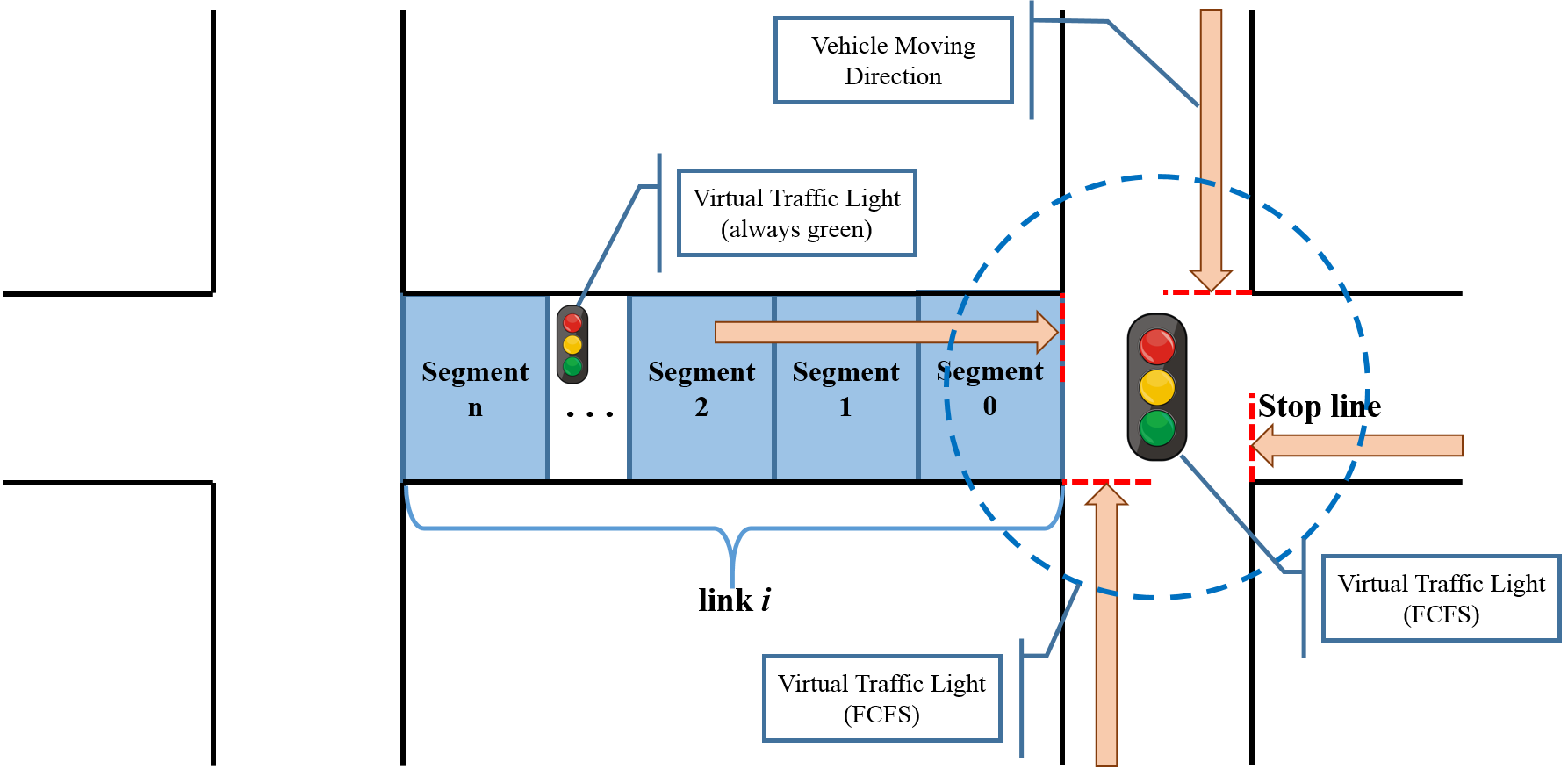}
	\caption{An AWSC Non-signalized Intersection}
	\label{fig_msm1}
\end{figure}

We make one more assumption. 
\begin{itemize}
	\item \textbf{A4:} Vehicles in each link move at the free flow speed whenever the distance between a vehicle in front and itself is more than the minimum separation distance.
\end{itemize}

\noindent Recall that in the previous subsection we propose a nonlinear flow rate function describing the intersection crossing model, where each driver needs to take an anticipation of the traffic light green time assignment into account. Since within a link there is no traffic light, which is equivant to the case that there has been an infinite number of continuous green time assignment, if we imagine there is some virtual traffic light inside each link, as shown in Figure 2, by that nonlinear flow rate model, we have that each driver intends to move with the highest allowable speed, only under the constraints of the minimum separation distance.  In this sense, Assumption \textbf{A4} and the nonlinear flow rate model described in the previous sub-section are consistent. 

With Assumption \textbf{A4}, we partition each link into segments such that the length of each segment is (roughly) equal to $v^*\times \lambda$, where $\lambda\in\mathbb{R}^+$ is the chosen period for each simulation step for each link connecting a non-signalized intersection. We assume that the sampling period $\Delta$ mentioned before is a multiple of $\lambda$. The reason why we choose a smaller simulation step period $\lambda$ is to improve the accuracy of the traffic model for each non-signalized intersection. Suppose a concerned link $i$ connecting to a non-signalized intersection consists of $n\in\mathbb{N}$ segments. We label them in an ascending order starting with the segment directly connecting with the non-signalized intersection with the segment ID of $0$, thus, the segment ID for the furthest segment from the non-signalized intersection is $n-1$. For each $j\in \{0,\cdots,n-1\}$ the link segment volume dynamic model $C_j(t)$ is the same as expression (1) shown before. But to distinguish it from the link volume dynamics related to signalized intersections, we use $t$ as the discrete time step variable. Thus, we have
\begin{equation}
	(\forall t\in\mathbb{N})\, C_i^j(t+1) = C_i^j(t) + d_i^j(t)-s_i^j(t),
	\label{eqn_seg_vehicle_1}
\end{equation}
where $d_i^j(t)$ and $s_i^j(t)$ are the number of entrance vehicles and exit vehicles of the segment $j$ in link $i$ at $t$, respectively. For segment 0 we define an indicator $T_i^0(t)\in\mathbb{N}$ such that $T_i^0(0):=+\infty$, and for all $t>0$ we have
\begin{equation} T_i^0(t+1):=\left\{\begin{array}{ll}\min\{T_i^0(t),t+1\} & \textrm{if $C_i^0(t+1)>0$}\\
		+\infty & \textrm{otherwise}\end{array}\right.
		\label{eqn_seg_indicator}
\end{equation}
The exit vehicle flow at segment 0 can be considered as a collection of continuous sequences, during each interval among the sequence segment volume is non-zero, separated by empty periods when there is no vehicle in the segment. For each $t\in\mathbb{N}$, there can be three possibilities: (1) $t$ is the leading interval of a continuous sequence, i.e., during $t-1$ there is no vehicle in segment 0. In this case $T_i^0(t):=t$, indicating when the current continuous sequence starts. (2) $t$ is part of a continuous sequence but is not the leading interval. Then clearly we have $T_i^0(t)=T_i^0(t-1)<t$, i.e., all intervals in the same continuous sequence get the leading interval as their indicators. (3) $t$ is an empty interval. In this case we have $C_i^0(t)=0$, and $T_i^0(t):=+\infty$. 

With such an indicator function for each link $i$ connecting to a non-signalized intersection $J\in\mathcal{J}_n$, we can define a virtual traffic light over the set of all links connecting to intersection $J$. Suppose each link $i\in\mathcal{L}_J$ has a unique ID $\sigma(i)\in\mathbb{N}$. We have the following traffic light green period assignment: for all $i\in\mathcal{L}_J$ and $t\in\mathbb{N}$, $\theta_i(t)=1$ if and only if the following condition holds:
\begin{equation}
\sigma(i)=\min\{\sigma(p)|p\in\mathcal{L}_J\wedge (\forall q\in\mathcal{F}_J)T_p^0(t)\leq T_q^0(t)\}
\label{eqn_seg_selection}
\end{equation}
that is, link $i$ gets a green light if and only if either none of the leading segments connecting to $J$ have vehicles and the ID of link $i$ is the smallest one (to break the tie), or there are links whose leading segments have vehicles waiting to cross $J$, and the vehicles in link $i$ reach the stop line at an earliest time and the ID of link $i$ is the smallest among all links whose vehicles reach the stop line at the same earliest time (to break the tie). Since the link ID is unique, there can be only one link getting the green time at each interval $t$. The difference between this intersection crossing model for a non-signalized intersection and a common FCFS model in our daily life is that in the latter case the size of each continuous sequence for segment 0 of each link is always 1, but in our model this size may be bigger than 1. The consequence is that, in our model it is likely that a vehicle in link $i$ comes earlier than a vehicle in link $i'$, but because the vehicle in $i'$ is part of a continuous sequence, whose leading interval is earlier than the arrival time of the vehicle in link $i$, the vehicle in link $i'$ actually gets the green light. In this sense we can see that our model imposes FCFS over continuous vehicle columns in relevant links instead of individual vehicles, i.e., once a column of vehicles gets the green light, it will continue having the green light until the column is over, i.e., in some interval $t$ the segment 0 of that link has no vehicles. The bottom line that forces us to consider a column-based FCFS, instead of an individual vehicle based FCFS is because our flow dynamic model cannot track individual vehicles, thus, each continuous column of vehicles are considered as one ``macro'' vehicle.

With the virtual traffic signals at these non-signalized intersections, we can adopt the flow dynamics constraints similar to the signalized intersections. For each segment $j$ in a link $i\in\mathcal{L}_J$ with $J\in\mathcal{J}_n$, the link outgoing shift $f_i^{j,j-1}(t)$ is determined as follows:
\begin{equation}
	\begin{aligned}
		& f_i^{j,j-1}(t) = \lfloor \min\{C_i^j(t),\hat{C}_i^{j-1}-C_i^{j-1}(t)+s_i^{j-1}(t),\\
		& \qquad  v_i^*d^*\lambda\}\rfloor
	\end{aligned}
		\label{eqn_seg_flow_dyna1}
\end{equation}
Note that there is no speed constraint here since all vehicles are assumed to run with the free speed within link $i$. The complication arises for the outgoing shift $f_{i,p}^0(t)$, where $p\in\mathcal{L}$ is the direct downstream link of segment 0 in link $i$. If $p$ is also part of a link connecting to a non-signalized intersection, we have
\begin{equation}
	\begin{aligned}
		& f_{i,p}^0(t) = \lfloor \min\{C_i^0(t),\hat{C}_p-C_p(t)+s_p(t),\\
		& \qquad  v_i^*d^*\lambda\}\rfloor
	\end{aligned}
	\label{eqn_seg_flow_dyna2}
\end{equation}
If $p$ is part of a signalized intersection, let $\Delta = m\lambda$, where $m\in\mathbb{N}$, we have the following expression:
\begin{equation}
	\begin{aligned}
		& f_{i,p}^0(t) = \lfloor \min\{C_i^0(t),\hat{C}_p-C_p(k)+s_p(k)-\\
		& \qquad \qquad\sum_{q=0}^{t\% m -1}f_{i,p}^0(t-1-q), v_i^*d^*\lambda\}\rfloor
	\end{aligned}
		\label{eqn_seg_flow_dyna3}
\end{equation}
where $t\% m$ denotes ``$t$ mod $m$''. In equation (15), because link $p$ directly connects a signalized intersection, its dynamic model uses the sampling period of $\Delta$, while the sampling period for link $i$ is $\lambda$. Thus, only after $m$ samples of link $i$, the dynamic model of link $j$ will be updated. For this reason we need to ensure that the vehicles coming from link $i$ within $\Delta$ will not overflow link $p$.

\section{A Formulation for Traffic Signal Control Problem for a Heterogeneous Traffic System}
\subsection{Objective function}

The total network-wise delay time within $N$ time intervals can be estimated as follows:
\begin{eqnarray}
\begin{aligned}
&\sum_{i\in \mathcal{L}}\sum_{k=1}^NC_i(k)\Big [1-\frac{\overline{v}_i(k)}{v_{i,max}}\Big ]\Delta \\
&= \sum_{i\in \mathcal{L}}\sum_{k=1}^N\Big [C_i(k)-\frac{L_i}{v_{i,max}}s_i(k)\Big ]\Delta
\label{eqn_delay_time}
\end{aligned}
\end{eqnarray}
where $\overline{v}_i(k)$ is the average speed of the link $i$ at $k$, which can be approximated by the ratio of the exit shift $s_i(k)$ and the average link density $C_i(k)/L_i$, where $L_i$ is the length of the link $i$. Here, we assume that vehicles in the link $i$ are uniformly distributed with identical speeds, i.e., the acceleration step is negligible. In reality, this assumption usually does not hold. But for the scheduling purpose, this is a sufficiently representative performance index. Since $s_i(k)=\sum_{j\in\mathcal{L}:(i,j)\in\cup_{J\in\mathcal{J}}\mathcal{F}_J}f_{ij}(k),$
we have our cost function as

\begin{equation*}
\min \sum_{i\in \mathcal{L}}\sum_{k=1}^N\Big [C_i(k)-\frac{L_i}{v_{i,max}}\sum_{j\in\mathcal{L}:(i,j)\in\cup_{J\in\mathcal{J}}\mathcal{F}_J}f_{ij}(k)\Big ]\Delta.
\label{eqn_cost_func}
\end{equation*}

\subsection{Conversions for logic and non-linear constraints}
\label{subsectionmld}
\subsubsection{Conversions for signalized intersections}
In the above description we have derived a linear cost function with a set of constraints describing traffic staging (\ref{eqn_RG_1})-(\ref{eqn_RG_3}), link volume dynamics (\ref{eqn_vehicle_1}) and exit shifts (\ref{eqn_speed_1})-(\ref{eqn_speed_5}), (\ref{eqn_fr_1}). Among these constraints, (\ref{eqn_RG_1}), (\ref{eqn_speed_1}), (\ref{eqn_speed_1.5}) and (\ref{eqn_speed_3}) are mixed logical constraints, which can be converted into mixed integer linear constraints by adopting a transformation strategy introduced in  \cite{Bemporad1999}, which is shown below.

Let $M$ be chosen to be sufficiently big, e.g., $M>\max_{i\in\mathcal{L}}\hat{C}_i$. Then (\ref{eqn_RG_1}) can be converted into
\begin{equation}
(\forall w\in\Omega_J)(\forall (i,j)\in h_J(w))\, f_{ij}(k)\leq M\theta_w(k).
\label{eqn_RG_1_conv}
\end{equation}
\begin{proposition}
	\textnormal{Replacing Condition (\ref{eqn_RG_1}) with Condition (\ref{eqn_RG_1_conv}) in the urban network traffic signal scheduling formulation leads to the same solution. \hfill $\Box$}
\end{proposition}
Proof: To see this conversion is valid, let $\theta_w(k)=0$, then $f_{ij}(k)\leq 0$. But since $f_{ij}(k)\geq 0$, we have $f_{ij}(k)=0$. If $\theta_w(k)=1$, it is trivially true that $f_{ij}(k)\leq M$ because by (\ref{eqn_fr_1}) we have $f_{ij}(k)\leq \max_{i\in\mathcal{L}}\hat{C}_i<M$.\hfill $\blacksquare$

Let $[0,r-1]:={0,\cdots,r-1}$. Conditions (\ref{eqn_speed_1}), (\ref{eqn_speed_1.5}) and (\ref{eqn_speed_3}) are equivalent to the following: for all $q\in[o,r]$,

\begin{subequations}
	\begin{align}
	&-(1-\theta_w(k-q-1))+ \delta_{ij}^{r-q}(k)\leq 0
	\label{eqn_speed_3_mixed_1}\\
	&(\forall p\in[0,q])-\theta_w(k-p)+ \delta_{ij}^{r-q}(k)\leq 0
	\label{eqn_speed_3_mixed_2}\\
	&(1-\theta_w(k-q-1))+\sum_{p=0}^q\theta_w(k-p) - \delta_{ij}^{r-q}(k)\leq q+1
	\label{eqn_speed_3_mixed_3}\\
	& l_{ij}(k)\leq \sum_{p=0}^r\delta_{ij}^p(k)l_{ij}^p + (r+1)(1-\theta_w(k+1))
	\label{eqn_speed_1_mixed_1}\\
	& l_{ij}(k)\geq \sum_{p=0}^r\delta_{ij}^p(k)l_{ij}^p - (r+1)(1-\theta_w(k+1))
	\label{eqn_speed_1_mixed_1.1}\\
	& l_{ij}(k)\leq \frac{1}{2}l_{ij}(k-1) + \theta_w(k+1)
	\label{eqn_speed_1.5_mixed_1.5}\\
	& l_{ij}(k)\geq \frac{1}{2}l_{ij}(k-1) - \theta_w(k+1)
	\label{eqn_speed_1.5_mixed_1.5.1}
	\end{align}
\end{subequations}

\begin{proposition}
	\textnormal{Replacing Conditions (\ref{eqn_speed_1}), (\ref{eqn_speed_1.5}) and (\ref{eqn_speed_3}) with Condition (\ref{eqn_speed_3_mixed_1})-(\ref{eqn_speed_1.5_mixed_1.5.1}) in the urban network traffic signal scheduling formulation leads to the same solution. \hfill $\Box$}
\end{proposition}
Proof: We can verify that Conditions (\ref{eqn_speed_3_mixed_1})-(\ref{eqn_speed_3_mixed_3}) are equivalent to Condition (\ref{eqn_speed_3}), Conditions (\ref{eqn_speed_1_mixed_1})-(\ref{eqn_speed_1_mixed_1.1}) are equivalent to Condition (\ref{eqn_speed_1}), and Conditions (\ref{eqn_speed_1.5_mixed_1.5})-(\ref{eqn_speed_1.5_mixed_1.5.1}) are equivalent to Condition (\ref{eqn_speed_1.5}).\hfill $\blacksquare$

{\color{black}
	
	By taking the cost function into account, which requires all shifts to be as large as possible, Condition (\ref{eqn_fr_1}) can be converted into the following mixed integer linear constraints:
	
	\begin{subequations}
		\begin{align}
		& f_{ij}(k)\leq \gamma_{ij}(k)C_i(k)
		\label{eqn_fr_1_mixed_logical_1}\\
		& f_{ij}(k)\leq l_{ij}(k)v_{i,max}d_{max}
		\label{eqn_fr_1_mixed_logical_2}\\
		& f_{ij}(k)\leq \hat{C}_j-C_j(k)
		\label{eqn_fr_1_mixed_logical_3}\\
		&\nonumber \big[f_{ij}(k) \geq \gamma_{ij}(k)C_i(k)\big] \vee \big[f_{ij}(k) \geq l_{ij}(k)v_{i,max}d_{max}\big]\\
		& \qquad \qquad \qquad \quad \quad \vee \big[f_{ij}(k) \geq \hat{C}_j-C_j(k)\big].
		\label {eqn_fr_1_mixed_logical_4}
		\end{align}
\end{subequations}}
Note that Condition (\ref{eqn_fr_1_mixed_logical_4}) is a logic contraint. However it can be removed from the formulation. Based on the objective function for this problem, the traffic flow rates will automatically choose the largest numbers they could reach. Thus this condition is redundant in the problem formulation. To simplify the formulation, we only take Conditions (\ref{eqn_fr_1_mixed_logical_1})-(\ref{eqn_fr_1_mixed_logical_3}) into consideration, which are all mixed integer linear constraints.
The merging condition (\ref{eqn_merging3}) is a linear inequality and does not require to change.
If we adopt the conservative merging condition shown in Eqn (\ref{eqn_merging}), it can be converted into the following mixed integer linear constraints:
\begin{subequations}
	\begin{align}
	& f^{(p)}_{ij}(k) \leq M\Gamma_{ij}(k)
	\label{eqn_mer_mixed_logical_1}\\
	& -f^{(p)}_{ij}(k) \leq -m\Gamma_{ij}(k)
	\label{eqn_mer_mixed_logical_2}\\
	& f^{(gw)}_{i'j}(k) + M\theta_i(k) + M \Gamma_{ij}(k) \leq 2M
	\label{eqn_mer_mixed_logical_3}
	\end{align}
\end{subequations}

\subsubsection{Conversions for non-signalized intersections}
For the formulation of the dynamics in non-signalized intersections, Eqn. (\ref{eqn_seg_indicator}) and (\ref{eqn_seg_selection}) are logic constraints and Eqn. (\ref{eqn_seg_flow_dyna1}) - (\ref{eqn_seg_flow_dyna3}) are non-linear constraints. For logic constraint (\ref{eqn_seg_indicator}), we first introduce a binary variable $\psi^1_i(t)$ which meet the following relationship,
\begin{eqnarray}
\begin{aligned}
& C_i^0(t)>0 \Rightarrow \psi^1_i(t) = 1\\
& C_i^0(t)=0 \Rightarrow \psi^1_i(t) = 0
\end{aligned}
\end{eqnarray}
which can be rewritten into the following mixed integer constraint set,
\begin{subequations}
	\begin{align}
	&-C_i^0(t)+ (M_C-1)\psi^1_i(t) \geq 0\label{eqn_seg_indicator_1}\\
	& C_i^0(t)- m_C\psi^1_i(t) \geq 0\label{eqn_seg_indicator_2}
	\end{align}
\end{subequations}
where $M_C = \max \{C_i^0(t)\}$ and $m_C = \min \{C_i^0(t)\}$. Note that Eqn. (\ref{eqn_seg_indicator_2}) is redundant and can be removed from the constraint set.
With $\psi^1_i(t)$ Eqn. (\ref{eqn_seg_indicator}) can be replaced by the following equation,
\begin{equation}
T_i^0(t+1) = M_T (1-\psi^1_i(t+1)) + \psi^1_i(t+1)\times \min\{T_i^0(t),t+1\}
\end{equation}
which can be replaced by some inequalities and logic constraints,
\begin{subequations}
	\begin{align}
		& T_i^0(t+1) \leq M_T (1-\psi^1_i(t+1)) + \psi^1_i(t+1)T_i^0(t)\label{eqn_seg_indicator_rep1}\\
		& T_i^0(t+1) \leq M_T (1-\psi^1_i(t+1)) + \psi^1_i(t+1)(t+1)\label{eqn_seg_indicator_rep2}\\
		&\nonumber [T_i^0(t+1) \geq M_T (1-\psi^1_i(t+1)) + \psi^1_i(t+1)T_i^0(t)]\\
		& \quad  \vee [T_i^0(t+1) \geq M_T (1-\psi^1_i(t+1)) + \psi^1_i(t+1)(t+1)]\label{eqn_seg_indicator_rep3}
	\end{align}
\end{subequations}
By introducing $\psi^2_i(t)$ and $\psi^3_i(t)$,
\begin{equation*}
\begin{aligned}
&[T_i^0(t+1) \geq M_T (1-\psi^1_i(t+1)) + \psi^1_i(t+1)T_i^0(t)]\\
&\qquad\leftrightarrow \psi^2_i(t+1) = 1\\
&[T_i^0(t+1) \geq M_T (1-\psi^1_i(t+1)) + \psi^1_i(t+1)(t+1)]\\
&\qquad\leftrightarrow \psi^3_i(t+1) = 1
\end{aligned}
\end{equation*}
Eqn. (\ref{eqn_seg_indicator_rep3}) is equivalent to the mixed integer constraint set below,
\begin{subequations}
	\begin{align}
	&\nonumber -m_T\psi^2_i(t+1) \leq T_i^0(t+1) - M_T (1-\psi^1_i(t+1))\\
	& \qquad\qquad\qquad\qquad - \psi^1_i(t+1)T_i^0(t) - m_T\\
	&\nonumber -(M_T+\epsilon)\psi^2_i(t+1) \leq - T_i^0(t+1) + M_T (1-\psi^1_i(t+1))\\
	&\qquad\qquad\qquad\qquad\qquad + \psi^1_i(t+1)T_i^0(t) - \epsilon\\
	&\nonumber -m_T\psi^3_i(t+1) \leq T_i^0(t+1) - M_T (1-\psi^1_i(t+1))\\
	& \qquad\qquad\qquad\qquad - \psi^1_i(t+1)(t+1) - m_T\\
	&\nonumber -(M_T+\epsilon)\psi^3_i(t+1) \leq - T_i^0(t+1) + M_T (1-\psi^1_i(t+1))\\
	&\qquad\qquad\qquad\qquad\qquad + \psi^1_i(t+1)(t+1) - \epsilon\\
	& \psi^2_i(t+1) + \psi^3_i(t+1) \geq 1
	\end{align}
\end{subequations}
Note in this constraint set, it includes some non-linear items such as $\psi^1_i(t+1)T_i^0(t)$ and $\psi^1_i(t+1)(t+1)$. Denote $\psi^1_i(t)T_i^0(t) = \psi^4_i(t)$ and $\psi^1_i(t)t = \psi^5_i(t)$, they can be transformed into the following linear constraints,
\begin{subequations}
	\begin{align}
	& \psi^1_i(t+1) \underline{T_i^0(t)} - \psi^4_i(t+1) \leq 0\\
	& \psi^4_i(t+1) - \psi^1_i(t+1) \overline{T_i^0(t)} \leq 0\\
	& T_i^0(t) + \psi^1_i(t+1) \overline{T_i^0(t)} - \psi^4_i(t+1) \leq \overline{T_i^0(t)}\\
	& \psi^4_i(t+1) - T_i^0(t) - \underline{T_i^0(t)} \psi^1_i(t+1) \leq - \underline{T_i^0(t)}\\
	& \psi^1_i(t+1) \underline{t} - \psi^5_i(t+1) \leq 0\\
	& \psi^5_i(t+1) - \psi^1_i(t+1) \overline{t} \leq 0\\
	& t + \overline{t} \psi^1_i(t+1) - \psi^5_i(t+1) \leq \overline{t}\\
	& \psi^5_i(t+1) - t - \underline{t} \psi^1_i(t+1) \leq - \underline{t}
	\end{align}
\end{subequations}

For Eqn. (\ref{eqn_seg_selection}) we have,
\begin{subequations}
	\begin{align}
	& (\forall q \in \mathcal{F}_J) T_p^0(t) \leq T_q^0(t)\\
	& [T_p^0(t) \geq T_{q_1}^0(t)] \vee [T_p^0(t) \geq T_{q_2}^0(t)] \vee \cdots
	\end{align}
\end{subequations}
By introducing some binary variables $\psi^*_i(t)\in\{\psi^6_i(t),\psi^7_i(t),\psi^8_i(t),\psi^{9}_i(t),\psi^{10}_i(t),\psi^{11}_i(t)\}$,
\begin{equation}
(\forall q \in \mathcal{F}_J) \psi^*_i(t) = 1 \leftrightarrow T_p^0(t) \geq T_{q}^0(t)
\end{equation}
Thus, they can be transformed into the mixed integer constraints,
\begin{subequations}
	\begin{align}
	 -m_T \psi^*_i(t) &\leq T_p^0(t) - T_{q}^0(t) - m_T\\
	 -(M_T - \epsilon) \psi^*_i(t) &\leq - T_p^0(t) + T_{q}^0(t) -\epsilon
	\end{align}
\end{subequations}
Since only one state can be true in each time interval, the summation of $\psi^6_i(t)$, $\psi^7_i(t)$, $\psi^8_i(t)$ and $\psi^9_i(t)$ is one.
\begin{equation}
\sum \psi^*_i(t) = 1
\end{equation}
Then the constraint to determine the virtual traffic signal is shown as follows,
\begin{equation}
(\forall q \in \mathcal{F}_J)[T_q^0(t) - T_p^0(t)=0] \leftrightarrow [\theta_i(t) = 1]
\end{equation}
which is equivalent to the inequalities shown below,
\begin{subequations}
	\begin{align}
	 -m_T\theta_i(t) & \leq T_q^0(t) - T_p^0(t) - m\\
	 - (M_T+\epsilon)\theta_i & \leq - T_q^0(t) + T_p^0(t) - \epsilon\\
	 -m_T\theta_i(t) & \leq T_p^0(t) - T_q^0(t) - m\\
	 - (M_T+\epsilon)\theta_i & \leq - T_p^0(t) + T_q^0(t) - \epsilon
	\end{align}
\end{subequations}

For Eqn. (\ref{eqn_seg_flow_dyna1}) - (\ref{eqn_seg_flow_dyna3}), please refer to the transformations for the flow dynamic constraints in signalized intersections.

{\color{black}When dealing with a large traffic network, solving the above MILP problem is certainly time consuming. To ensure a real-time solution to this problem, we present a distributed matheuristic approach in our previous work \cite{Zhang2015}.}

\section{Simulation Results}
The numerical simulations for this proposed model and traffic signal control problem are done in MATLAB with Gurobi solver. Alternatively, some real-time traffic simulations are done in VISSIM, a microscopic traffic simulator. 

%A comparison between the results generated by the fixed cycle traffic scheduling strategy in VISSIM and the adaptive traffic signal scheduling strategy based on the proposed flow dynamic model are shown in Fig. \ref{fig_vm_res_3}. The light blue bar is the network-wise average delay for the fixed time strategy with reversing time set to 30s and the red bar is the the network-wise average delay for the fixed time strategy with 60s reversing time. The dark blue bar shows the delay time after adopting our proposed system model and strategy with the sampling time of 15s. Results show that our proposed approach can significantly reduce (more than 60\% in average) the network-wise time delay.
%
%\begin{figure}[!ht]
%	\centering
%	\includegraphics[width = \linewidth]{vissim_results_3.png}
%	\caption{A comparison between fixed cycle strategy in VISSIM and the adaptive strategy with the proposed traffic model}
%	\label{fig_vm_res_3}
%\end{figure}

\subsection{Simulation-based Model Validation in VISSIM}
We use VISSIM to build an urban traffic network and validate the proposed non-signalized intersection model. 
\subsubsection{The VISSIM simulation procedure} 
For a 4-arm non-signalized intersection, there are two groups of virtual traffic signals at each intersection, i.e., the signal group for horizontal links and the signal group for vertical links, both with bi-directional traffic flows.
The VISSIM simulation procedure is shown as follows. First, we set the virtual traffic signal at the non-signalized intersection to be RED for all directions and wait for the first vehicle's coming. We set the traffic signal to be GREEN for one signal group if we detect the vehicle's coming from current direction. The traffic signal will be kept in the next $t_{acp}$ seconds and at the end of the GREEN signal, we will detect if the link segment related to this GREEN signal is empty. If so, we give GREEN signal to the other signal group; if not, we keep the GREEN signal for this direction until the link segment is cleared.

%\noindent \textbf{Step 1: Setup}
%\begin{enumerate}
%	\item All detectors in the horizontal direction(left to right \& right to left) are treated as 1, no matter how many lanes there are. At each simulation step, if detector=1, save the trigger time for this direction, so as to the vertical direction.
%\end{enumerate}
%
%\noindent\textbf{Step 2: Initialization}
%\begin{enumerate}
%	\item All virtual traffic signals are set to RED. 
%	\item Start to inject vehicles, a signal turns GREEN when the corresponding direction gets cars first. 
%	\item If both directions get cars simultaneously, the horizontal direction turns GREEN to break the tie.
%\end{enumerate}
%
%\noindent\textbf{Step 3: Iterative updating of virtual signal values}
%\begin{enumerate}
%	\item If the virtual signal is GREEN, start to check whether the queue is cleared by checking whether a car has arrived within the next 2 seconds. If \textit{true}, the queue is not cleared and we continue to give GREEN for this direction. If there are no car coming within the next 2 seconds, proceed to turn GREEN to RED.
%\end{enumerate}

\subsubsection{The validation of non-signalized intersection model}
To validate the proposed model, we compare the results from the VISSIM simulator with the results predicted by our non-signalized intersection model.
Two cases are generated for this validation, i.e., the low volume and high volume scenarios, respectively. 
\begin{figure}[!ht]
	\centering
	\includegraphics[width=0.9\linewidth]{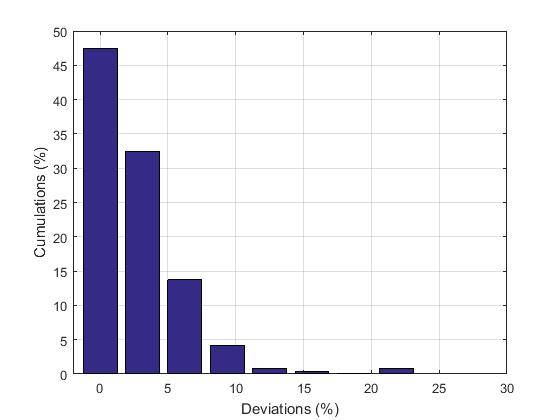}
	\caption{The distribution for the deviations between the VISSIM simulation data and the analytical data for low traffic volume case}
	\label{fig:480}
\end{figure}
\begin{figure}[!ht]
	\centering
	\includegraphics[width=0.9\linewidth]{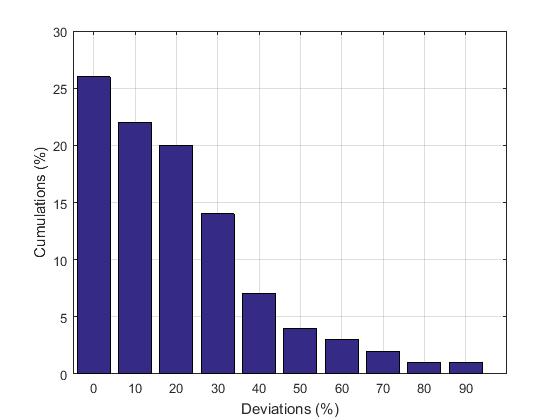}
	\caption{The distribution for the deviations between the VISSIM simulation data and the analytical data for high traffic volume case}
	\label{fig:1200}
\end{figure}
The discrete time interval is set to be 15s and the whole simulation period is 3600s. There are 240 sampling points included in each case. For each sampling point, we obtain the incoming flow for each link connected to the non-signalized intersection and outgoing flow at this intersection in VISSIM. In the mean time, we pass the same incoming flow from VISSIM to MATLAB and calculate the outgoing flow by the proposed non-signalized intersection model. The deviations between these two groups of results are shown in Fig. \ref{fig:480} and Fig. \ref{fig:1200}.
Fig. \ref{fig:480} shows the deviations in the low traffic volume case. We can see that the maximal deviation in this case is no more than $25\%$ and about $80\%$ of data are under $5\%$ deviation.
$97.92\%$ of data are under $10\%$ deviation, which shows the effectiveness of the proposed non-signalized flow dynamics model. 
Fig. \ref{fig:1200} shows the deviations in the high traffic volume case. In this case, there does exist some data with very large deviations because of the approximation nature of our proposed model, as mentioned in Sub-section II.D, where we impose FCFS over columns of vehicles instead of each single one.  The figure shows that about $50\%$ of the data are under $10\%$ deviation and about $85\%$ of the data are under $30\%$ deviation.
As a conclusion, the proposed non-signalized intersection model can provide a relatively precise description for the flow dynamics with low traffic volumes. How to develop a more accurate flow dynamic model for a non-signalized intersection with high traffic volumes remains an interesting research topic, which hopefully could be addressed  in our further work.

\subsection{Simplified Traffic Network Profile}
\begin{figure}[!ht]
	\centering
	\includegraphics[width=0.9\linewidth]{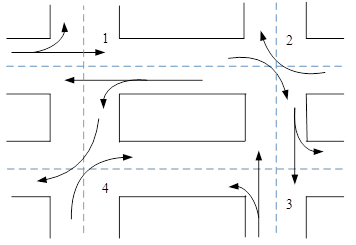}
	\caption{4-stage system}
	\label{fig_4stage}
\end{figure}

Fig \ref{fig_4stage} shows four stages of a signalized intersection, which is currently utilized in the urban traffic network of Singapore. Different stages include a set of different traffic flows. 
In the fixed cycle traffic scheduling strategy, the four stages are repeated based on the order $1-2-3-4$. In our model, each stage is assigned with a period in a real-time manner by the network controller with no cycle based reappearance. We use four integers, ``1", ``2", ``3" and ``4", i.e., $w\in\Omega_J = \{1,2,3,4\}$, to represent four stages of one intersection shown in Fig \ref{fig_4stage}. 

\begin{figure}[!ht]
	\centering
	\includegraphics[width=0.8\linewidth]{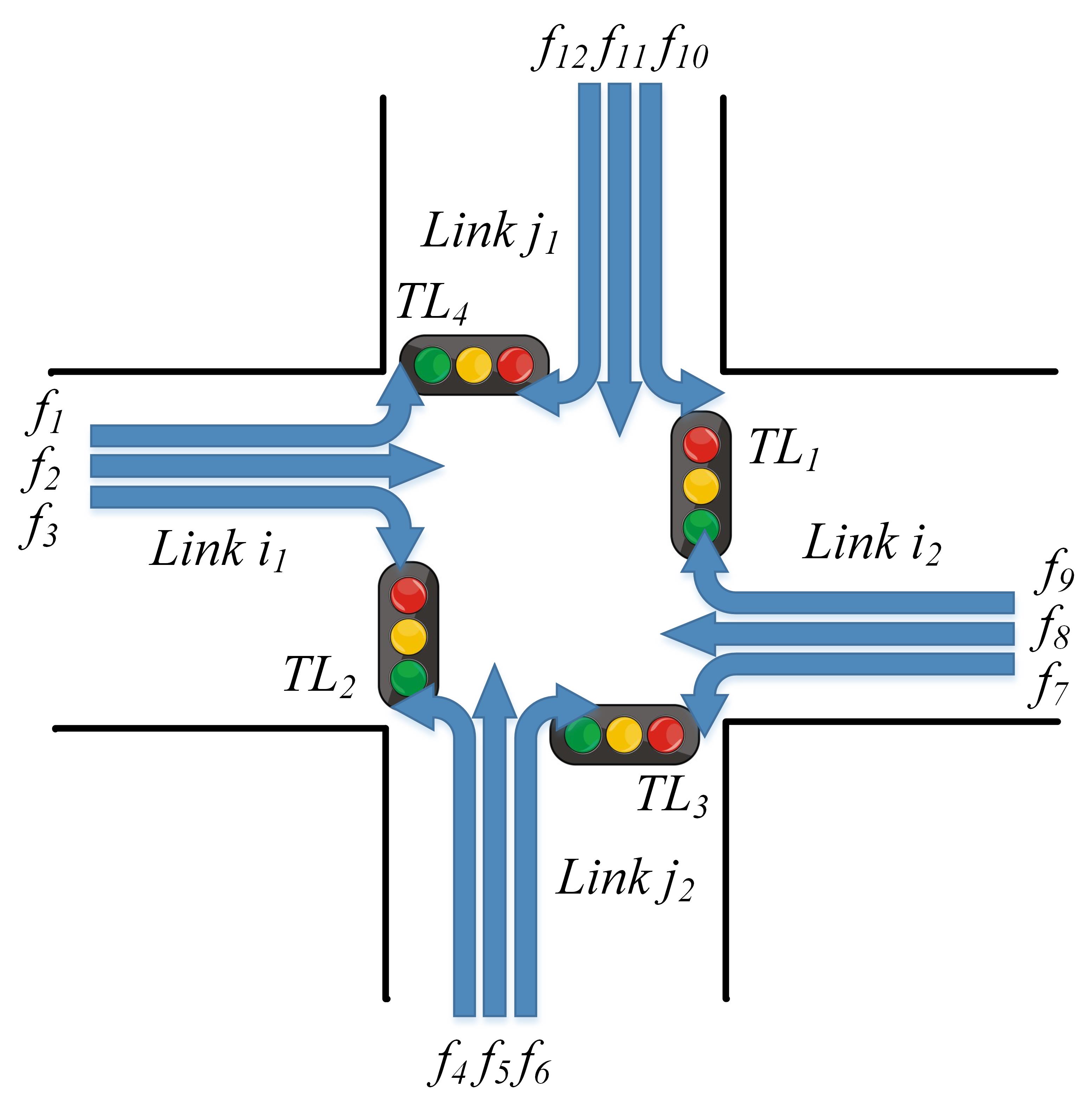}
	\caption{An example for bi-directional intersection}
	\label{fig_msm}
\end{figure}

For this bi-directional four-stage system, each link contains 3 traffic flows, i.e., straight-forward, left-turning and right-turning , thus each intersection contains a total of 12 traffic flows. For the intersection $J\in\mathcal{J}$ which connected with links $i_i$, $i_2$, $j_1$ and $j_2$, the flows are denoted as $f_{ij}$, where $i\in\{i_1, i_2, j_1, j_2\}$ and $j\in\{i_1, i_2, j_1, j_2\}\backslash\{i\}$, as shown in Fig. \ref{fig_msm}.

\subsection{A Comparison of Fully Controlled Traffic Network, Partially Controlled Traffic Network and Fully Uncontrolled Traffic Network}
We introduce a numerical comparison generated in MATLAB on the system performance among the fully controlled traffic network, partially controlled traffic network and fully uncontrolled traffic network, i.e., a traffic network with all intersections signalized, part of the intersections signalized and all intersections non-signalized. 
Nine cases with different traffic densities are created in this case study as shown in Table \ref{tab_case_study}.
\begin{table}[!ht]
	\centering
	\caption{A Comparison of Fully Signalized Traffic Network and Partially Signalized Traffic Network}
	\label{tab_case_study}
	\begin{tabular}{|c|p{0.18\linewidth}|p{0.18\linewidth}|p{0.2\linewidth}|}
		\hline
		
		& Fully Controllable & Partially Controllable & Fully uncontrollable\\
		
		\hline
		
		Low Traffic & Case 1 & Case 2 & Case 3\\
		
		\hline
		
		Medium Traffic & Case 4 & Case 5 & Case 6\\
		
		\hline
		
		High Traffic & Case 7 & Case 8 & Case 9\\
		
		\hline
		
	\end{tabular}
\end{table}
For the \textit{Low Traffic} cases, we mean that the traffic volumes for each direction are lower than the boundaries of the gap-acceptance model; for the \textit{Medium Traffic} case, the traffic volumes on some links may be higher than the gap-acceptance model and for the \textit{High Traffic} case, all the traffic volumes are out of the boundaries of the gap-acceptance model, which may lead to high congestions at some intersections.
All these cases are tested based on a 4-to-4 traffic grid and a total of 16 intersections are involved in this system.
For the fully controlled cases, we put 16 traffic signals in this system and for the partially controlled cases 9 of the 16 intersections are signalized. The sampling time is set to be 15s and the total simulation time period is 180s.
The network-wise time delays are calculated and shown in Table \ref{tab_case_results}.
\begin{table}[!ht]
	\centering
	\caption{A comparison of the network-wise traffic delays for the nine cases}
	\label{tab_case_results}
	\begin{tabular}{c|c|c|c}
		\hline
		
		Case No. & Total vehicles & Total delays (sec) & Average delays (sec) \\
		
		\hline 
		
		Case 1 & 320 & 0 & 0\\
		
		\hline
		
		Case 2 & 320 & 4860 & 486\\
		
		\hline
		
		Case 3 & 320 & 8640 & 864\\
		
		\hline
		
		Case 4 & 560 & 11700 & 1170\\
		
		\hline
		
		Case 5 & 560 & 11520 & 1152\\
		
		\hline
		
		Case 6 & 560 & 9360 & 936\\
		
		\hline
		
		Case 7 & 800 & 17640 & 1764\\
		
		\hline
		
		Case 8 & 800 & 15360 & 1536\\
		
		\hline
		
		Case 9 & 800 & 13200 & 1320\\
		
		\hline
		
	\end{tabular}
\end{table}

\begin{figure}[t]
	\centering
	\includegraphics[width = 0.95\linewidth]{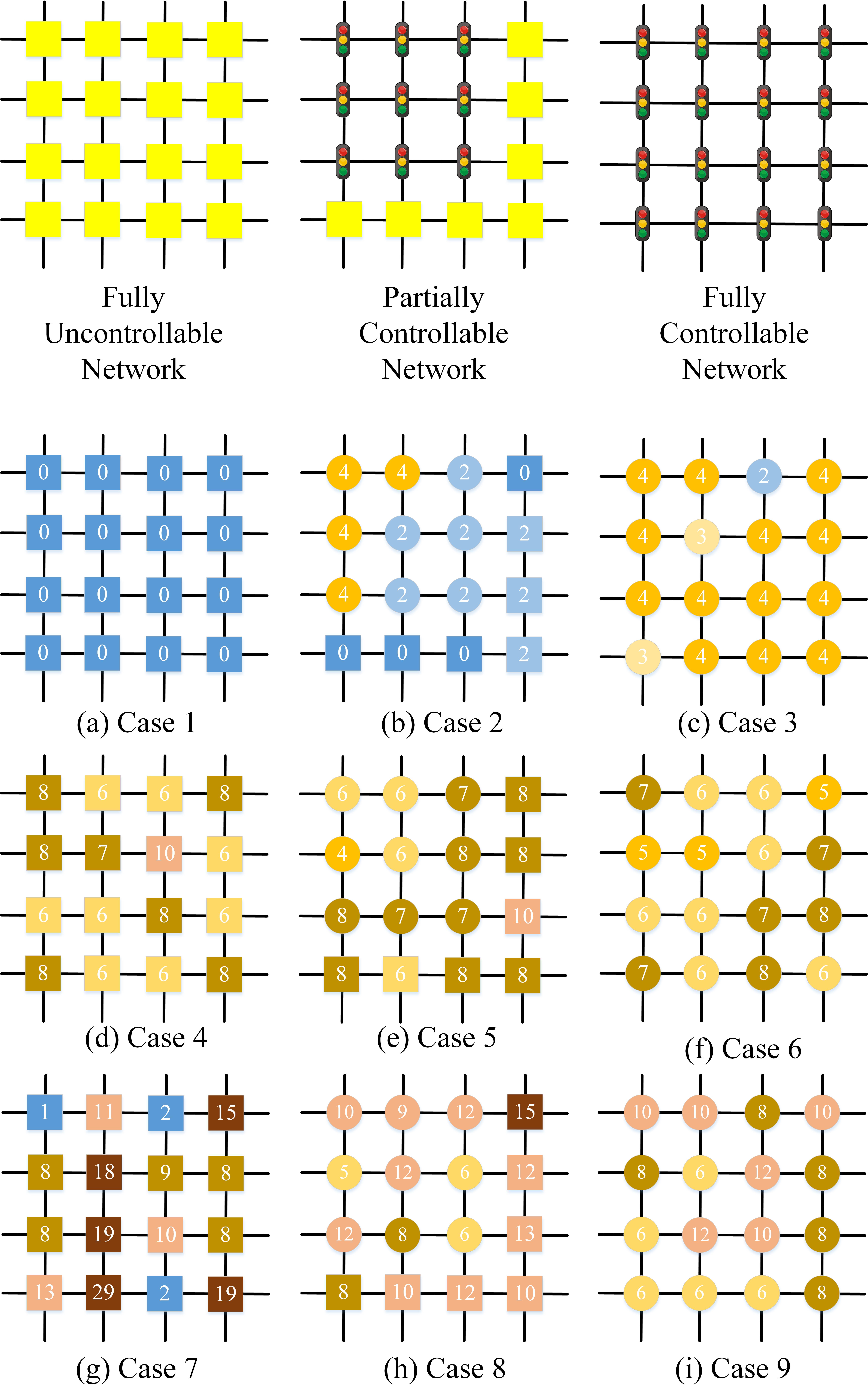}
	\caption{Number of vehicles waiting at the intersections in the 9 case studies}
	\label{fig_cases}
\end{figure}

For Case 1, we assume that the vehicles on all the links will not touch the threshold for the gap-acceptance model, which means that all the vehicles from different directions can pass the non-signalized intersections within one time interval unless the downstream links are blocked. Thus the network-wise traffic delay for Case 1 is 0. For Case 2, we put 9 signalized intersections. It means that for those 9 intersections only one stage can be activated for each time interval, which will block the traffic from the other direction. In this case, the network-wise traffic delay is 4860s. If we signalize all the intersections as shown in Case 3, the traffic delay will keep increasing because the freedom of flow transmission is reduced with the increasing of traffic signals.
As a conclusion, for the low traffic cases, traffic signals may introduce unnecessary traffic delays to the whole system.
For the \textit{Medium Traffic} cases and \textit{High Traffic} cases, the introducing of traffic signals does reduce the traffic delays in the system. In some scenarios such as Case 5, the signalized intersections do not have much more improvements. However, when the traffic congestions are heavy, the signalized intersections will help to significantly reduce the delay time, as shown in Case 8 and Case 9.

In Fig. \ref{fig_cases}, we demonstrate the vehicles delayed at the signalized or non-signalized intersections. The circles denote the signalized intersections and the rectangles denote the non-signalized intersections. The color changes from dark blue to dark brown, which denotes the level of congestions from low to high. Numbers in the circles or rectangles are the number of vehicles delayed at the intersections from all the directions.

\subsection{VISSIM Simulation in Singapore Traffic Networks}
\begin{figure}[!ht]
	\centering
	\includegraphics[width=0.85\linewidth]{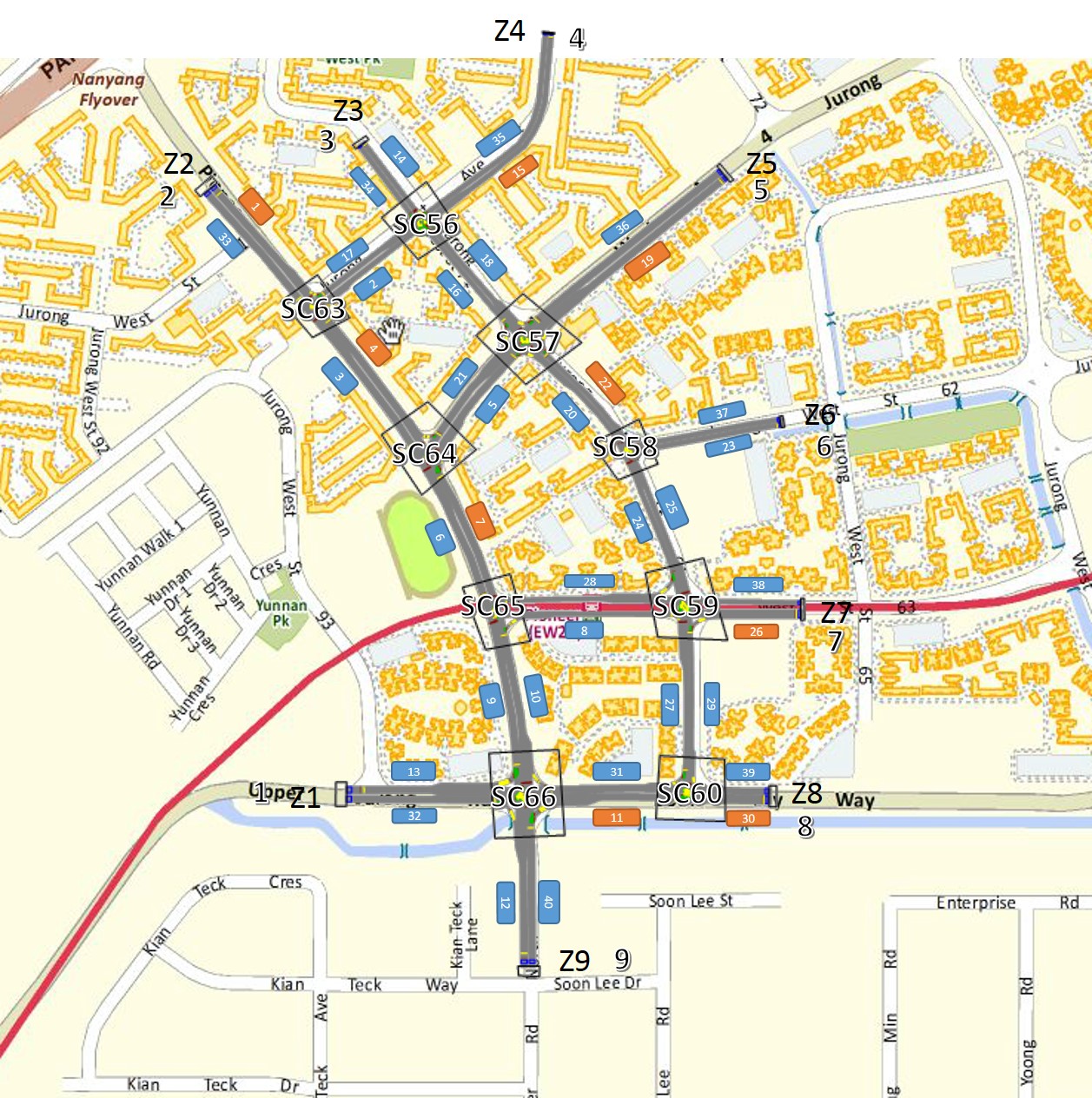}
	\caption{VISSIM Simulation Platform built in Jurong West in Singapore (The map file is downloaded from https://www.onemap.sg/index.html.)}
	\label{fig:9}
\end{figure}

A simulation platform is developed in VISSIM based on Jurong West Area in Singapore which is shown in Fig. \ref{fig:9}, which consists of 9 intersections and 40 links. 
The VISSIM simulation platform is connected to MATLAB via COM interface and the traffic signal controller is developed in MATLAB with Gurobi optimization solver.
We test this case study with fully uncontrollable network, partially controllable network and fully controllable network with low traffic, medium traffic and high traffic volumes as well.

%Note that for this case study, the computational complexity is acceptable since the number of involved intersections is not so big. To deal with large-scale traffic network, a distributed computational structure \cite{Zhang2015} with meta-heuristics \cite{Gao2017} \cite{Gao2016} can be applied to reduce the complexity and to meet the requirement of real-time processing. 

%A simulation platform is developed in VISSIM based on Jurong Area in West Singapore which is shown in Fig. \ref{fig:9}.
%This area includes 66 intersections, including signalized intersections and non-signalized intersections.
%We test our traffic signal scheduling algorithm for a heterogeneous urban traffic system with MATLAB and VISSIM. The traffic signal controller is developed in MATLAB and the control signal is transferred to VISSIM via COM-interface. 
%
%\begin{figure}[!ht]
%	\centering
%%	\includegraphics[width=0.95\linewidth]{9}
%	\includegraphics[width=0.95\linewidth]{Picture1}
%%	\includegraphics[width=0.95\linewidth]{Picture3}
%	\caption{VISSIM Simulation Platform}
%	\label{fig:9}
%\end{figure}
%

\section{Conclusions}
In this paper, a novel model for the urban traffic system is proposed to describe the heterogeneous traffic system with signalized intersections and non-signalized intersections. The proposed model is validated based on the simulations in VISSIM. Moreover, an urban traffic signal control problem formulation for the traffic network with signalized and non-signalized intersections is proposed as a mixed integer programming problem. Some comparisons among the fully controlled, partially controlled and fully uncontrolled traffic system are provided based on the numerical analyses and we show the potential applications in the traffic system design, which leaves the door open for developing a systematic planning approach on deciding what traffic junctions require signal control to ensure a good traffic control performance, thus, have a great social and economic potentials, considering that it is rather expensive to have signal control in an urban area.

\begin{IEEEbiography}[{\includegraphics[width=1in,height=1.25in,clip,keepaspectratio]{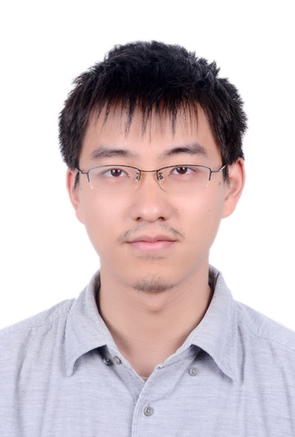}}]{Yicheng Zhang (S'15)}
	Yicheng Zhang received his Master degree from University of Science and Technology of China in 2014. He is currently a PhD candidate in School of Electrical and Electronic Engineering at Nanyang Technological University, Singapore. His current research interests include optimization method, model-based fault diagnosis, smart grid and intelligent transportation systems.
\end{IEEEbiography}

\begin{IEEEbiography}[{\includegraphics[width=1in,height=1.25in,clip,keepaspectratio]{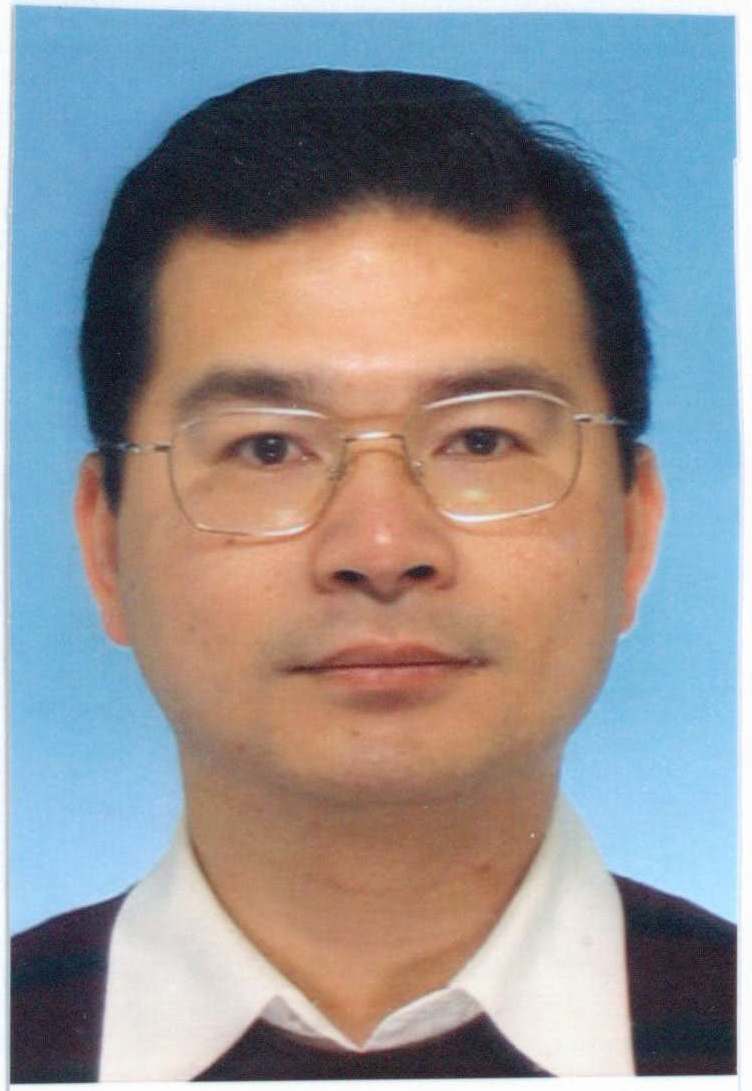}}]{Rong Su (M'11-SM'14)}
	Dr Su Rong obtained his Bachelor of Engineering degree from University of Science and Technology of China, and Master of Applied Science degree and PhD degree from University of Toronto. Since then he was affiliated with University of Waterloo and Technical University of Eindhoven before he joined Nanyang Technological University in 2010. Dr Su's research interests include discrete-event system theory, model-based fault diagnosis, operation planning and scheduling and multi-agent systems with applications in flexible manufacturing, intelligent transportation, human-robot interface, power management and green building. So far he has been involved in several projects sponsored by Singapore National Research Foundation (NRF), Singapore Agency of Science, Technology and Research (A*STAR), Singapore Ministry of Education (MoE), Singapore Civil Aviation Authority (CAAS) and Singapore Economic Development Board (EDB). Dr Su is a senior member of IEEE, and an associate editor for Transactions of the Institute of Measurement and Control and Journal of Control and Decision. He is also the Chair of the Technical Committee on Smart Cities in the IEEE Control Systems Society.
\end{IEEEbiography}

\begin{IEEEbiography}[{\includegraphics[width=1in,height=1.25in,clip,keepaspectratio]{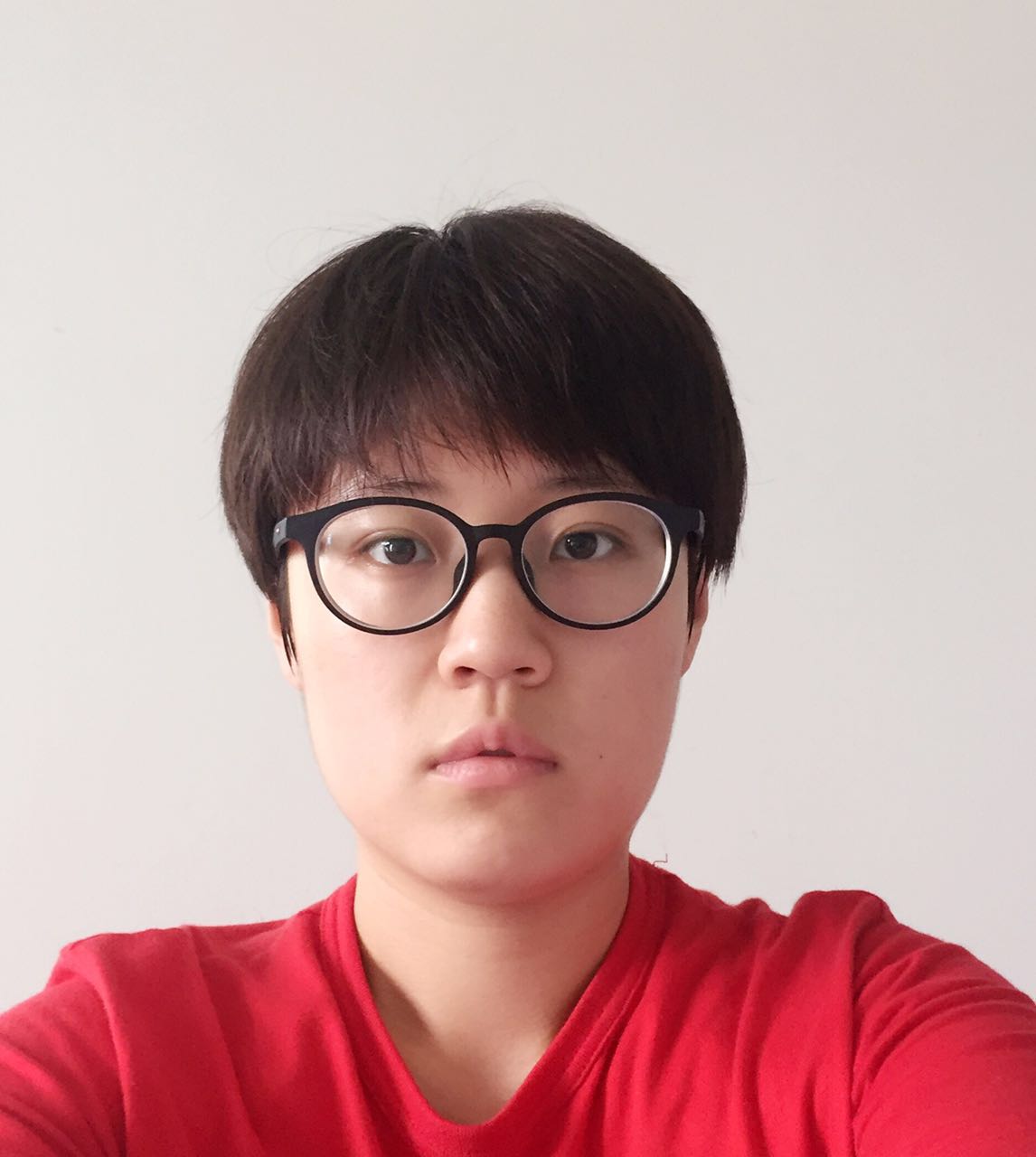}}]{Yi Zhang (S'17)}
	Yi Zhang received her Bachelor degree from China in Shandong University in 2014. She is currently a PhD candidate at Nanyang Technological University, Singapore. Her current research interests include model-based traffic signal scheduling and multi-directional pedestrian-flow simulation.
\end{IEEEbiography}

\begin{IEEEbiography}[{\includegraphics[width=1in,height=1.25in,clip,keepaspectratio]{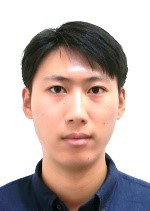}}]{Chunyang Sun}
	Chunyang Sun received B.Eng. (EEE) from Nanyang Technological University (NTU), Singapore, 2016. He was a software engineer in Innogrity Pte Ltd (Singapore) from Sep. 2012 to Aug. 2014. Since January 2016, he was a project officer in the School of Electrical and Electronic Engineering (EEE), NTU, Singapore. His research interests include traffic simulation and software development.
\end{IEEEbiography}

\vfill

\end{document}